\documentclass[a4paper]{article}
\usepackage[T1]{fontenc}
\usepackage{url}
\usepackage[dvipdfmx]{graphicx}
\usepackage{latexsym}
\usepackage[numbers]{natbib}
\usepackage{rotating}
\usepackage{amsmath,amssymb,amsfonts}
\usepackage{algorithm,algorithmic}
\usepackage{textcomp}
\usepackage{xcolor}

\begin{document}
%
\title{Performance evaluation of accelerated complex multiple-precision LU decomposition}
%
%
\author{Tomonori Kouya\\Shizuoka Institute of Science and Technology\\OrcidID:0000-0003-0178-5519}
%
%
%
\maketitle              
\begin{abstract}
The direct method is one of the most important algorithms for solving linear systems of equations, with LU decomposition comprising a significant portion of its computation time. This study explores strategies to accelerate complex LU decomposition using multiple-precision floating-point arithmetic of the multiple-component type. Specifically, we explore the potential efficiency gains using a combination of SIMDization and the 3M method for complex matrix multiplication. Our benchmark tests compare this approach with the direct method implementation in MPLAPACK, focusing on computation time and numerical errors.
\end{abstract}
%
\section{Introduction}

Linear systems of equations are frequently encountered in various scientific computations, therefore various solvers have been developed to address the computation environment and mathematical properties of these equations. The direct method, including LU decomposition, is considered standard because it can solve equations with a limited number of algebraic operations and maintain numerical stability through the use of a pivoting strategy. The HPL benchmark test\cite{hpl_linpack} employs the direct method to select the top 500 supercomputers concurrently.

In recent years, there has been much research on mixed-precision computing, where the number of precision digits (the length of mantissa) of the floating-point (FP) numbers used can be tailored to the requirements of the problem. This approach allows for improving computation performance by combining various kinds of hardware-oriented IEEE754 FP arithmetic, such as binary16 (11-bit precision), binary32 (24-bit), and binary64 (53-bit). Nevertheless, inadequate accuracy in calculation results would pose a significant disadvantage. There is an increasing demand for ``reproducible computing" to minimize environmental dependence of the calculation results. For good-condition problems, speed-up should be pursued by using fast binary16, 32, and 64 FP arithmetic, which have shorter calculation digits. Conversely, for bad-condition problems, multiple-precision FP operations exceeding binary64 should be employed to ensure accuracy and reproducibility. In particular, for stable solutions of nonlinear problems based on linear computation, it is necessary to minimize rounding errors in the intermediate results. This, utilizing a standard multiple-precision FP arithmetic library known for its user-friendly interface and high performance, such as QD\cite{qd} and MPFR\cite{mpfr} becomes essential. MPLAPACK/MPBLAS\cite{mplapack}, serving as a de-facto standard multiple-precision linear computation environment incorporating these libraries, becomes essential. We have developed three high-performance multiple-precision real basic linear calculations of double-double (DD, 106 bits), triple-double (TD, 159 bits), and quadruple-double (QD, 212 bits) using AVX2\cite{kouya_iccsa2021}. The results have been compiled into BNCmatmul\cite{bncmatmul} and released as open-source software.

Based on these results, this paper discusses the performance evaluation of complex multiplicative precision LU decomposition, particularly focusing on multi-component DD, TD, and QD precision, while also employing the 3M method, a type of complex multiplication method. Benchmark tests compared the acceleration of complex arithmetic using AVX2, the LU decomposition using the Strassen algorithm and the Ozaki scheme, and parallelization using OpenMP.

\section{Acceleration of complex basic linear computation}

As mentioned above, we have already implemented real basic linear computation and LU decomposition with SIMD operations such as AVX2 and a library of multiple-precision linear computation supporting parallelization with OpenMP. Additionally, we have developed multi-component fixed-precision operations supporting DD, TD, and QD precision, which are faster than MPLAPACK/MPBLAS, and arbitrary precision floating-point operations based on MPFR. Building upon this foundation, we have implemented complex BLAS\cite{kouya_complex_lu_mpfr}. It is worth noting that while we have published a collection of results on arbitrary-precision complex LU decomposition\cite{kouya_complex_lu_mpfr}, we have yet to parallelize it. In the following, we explain the strategies for improving the performance of complex linear computation, which serves as the basis for complex LU decomposition.

\subsection{The 3M method for complex linear computation}

For any $a, b\in\mathbb{C}$, standard complex multiplication $ab\in\mathbb{C}$ is calculated using the 4M (multiplication) method, which involves four real multiplications:
\begin{equation}
\begin{split}
    ab &:= (\mathrm{Re}(a)\mathrm{Re}(b) - \mathrm{Im}(a) \mathrm{Im}(b)) + \left\{ \mathrm{Im}(a)\mathrm{Re}(b) \right.\\
    &\left. + \mathrm{Re}(a)\mathrm{Im}(b)\right\} \cdot \mathrm{i},\ 
\end{split}\label{eqn:4m}
\end{equation}
In contrast, the 3M method is executed as follows:
\begin{equation}
    \begin{split}
        t_1 &:= \mathrm{Re}(a)\mathrm{Re}(b),\ t_2 := \mathrm{Im}(a) \mathrm{Im}(b) \in\mathbb{R}\\
        ab &:= (t_1 - t_2) + \left\{ (\mathrm{Re}(a) + \mathrm{Im}(a))\right. \\
        & \left.(\mathrm{Re}(b) + \mathrm{Im}(b)) - t_1 - t_2 \right\}\cdot \mathrm{i} \in\mathbb{C},
    \end{split}\label{eqn:3m}
\end{equation}
The 3M method retains the same computation in the real part as the 4M method but reduces the number of real multiplications in the imaginary part, hence its name. We cannot ascertain the efficiency of the 3M method without benchmark testing of its implementations because it includes three additional additions compared to the 4M method. Nevertheless, The 3M method has proven to be valuable for tasks such as matrix multiplications\cite{cmatmul_3m4m} and arbitrary-precision complex arithmetic\cite{mpc}\cite{mpsolve2000}, particularly in scenarios with extensive computation times for multiplication.

We implement multiple-precision matrix multiplication based on the 3M method. For any 
 $A\in\mathbb{C}^{m\times n}$ and $B\in\mathbb{C}^{n\times l}$:  
\begin{equation}
    \begin{split}
        T_1 &:= \mathrm{Re}(A)\mathrm{Re}(B),\ T_2 := \mathrm{Im}(A)\mathrm{Im}(B)  \in\mathbb{R}^{m\times l} \\
        AB &:= (T_1 - T_2) + \left\{ (\mathrm{Re}(A) + \mathrm{Im}(A)) \right.\\
        & \left.(\mathrm{Re}(B) + \mathrm{Im}(B)) - T_1 - T_2 \right\}\cdot\mathrm{i}\in\mathbb{C}^{m\times l},
    \end{split}\label{eqn:3m_cemm}
\end{equation}
This approach allows for complex matrix multiplication based on existing high-performance real matrix multiplication techniques. While the efficiency of this implementation remains uncertain due to the additional three additions and subtractions of matrices, benchmark tests are necessary to evaluate the performance of our multiple-precision complex matrix multiplication.

\subsection{Acceleration using AVX2 SIMDization}

\begin{figure}[htb]
    \begin{center}
    \includegraphics[width=.9\textwidth]{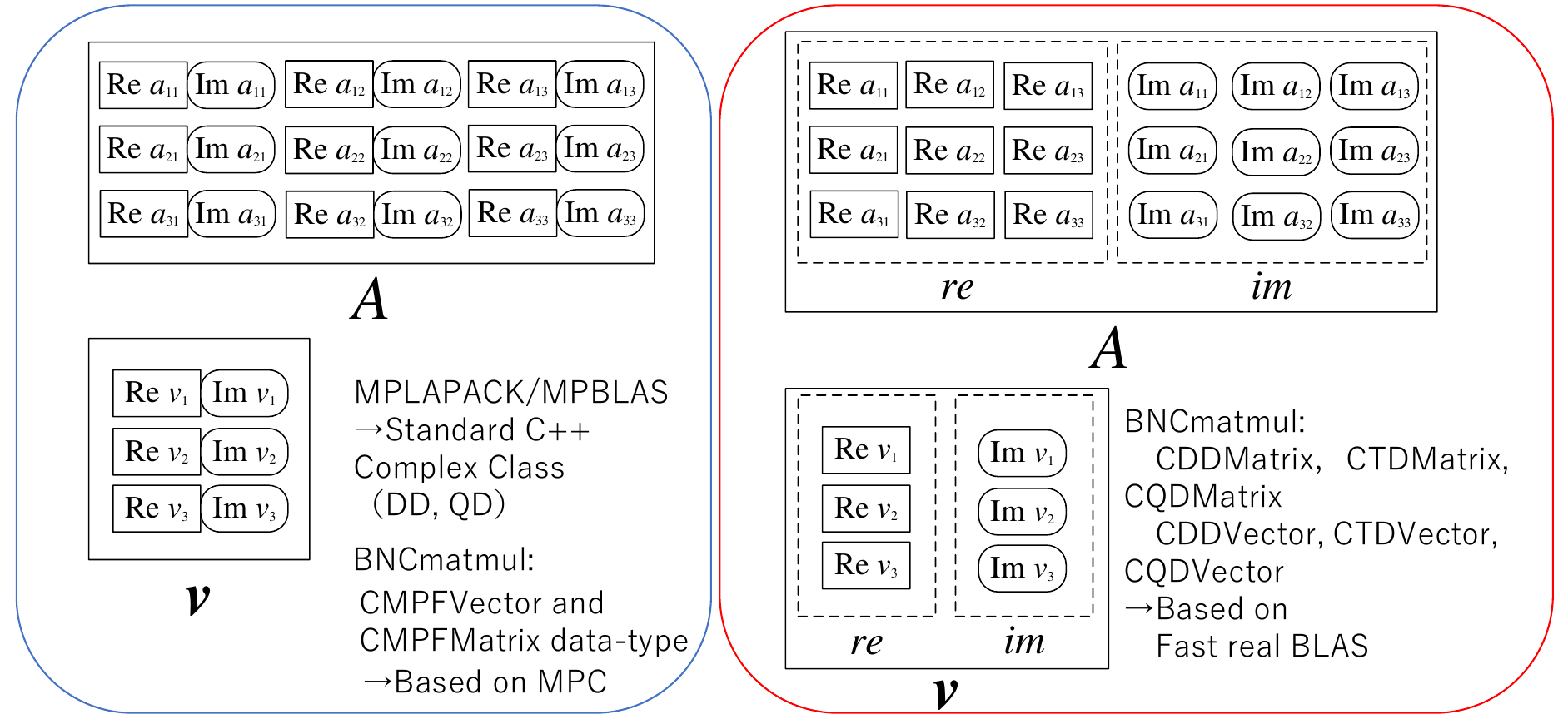}
    \caption{Data structure for multiple-precision complex vectors $\mathbf{v}\in\mathbb{C}^3$ and matrices $A\in\mathbb{C}^{3\times 3}$}\label{fig:structure_complex_matrix_vector}
    \end{center}
\end{figure}

As described above, our complex basic linear computation, which includes the 3M method is constructed upon real multiple-precision basic linear subprograms (BLAS) developed with SIMDization such as AVX2 in our library. Complex vectors and matrices are allocated as separate arrays for real and imaginary parts in multi-component-type multiple-precision FP data types. On the other hand, for arbitrary-precision FP arithmetic, we adopt the same data type (a set of real and imaginary parts) as MPC\cite{mpc} for a single element. MPC, based on MPFR, employs the 3M method for element-wise complex multiplication, ensuring high speed simply by using it. The data structure for these complex vectors and matrices is illustrated in \figurename\ref{fig:structure_complex_matrix_vector}.

This approach simplifies the implementation of fast complex matrix multiplication based on the formula in (\ref{eqn:3m_cemm}) for multi-component-type multiple-precision calculations. We have already implemented four matrix multiplication algorithms: simple matrix multiplication using triple loops, block matrix multiplication, Strassen matrix multiplication (and the Winograd algorithm), and the Ozaki scheme. Complex matrix multiplication is achieved by combining them.

\begin{figure*}[htb]
    \begin{center}
    \includegraphics[width=.9\textwidth]{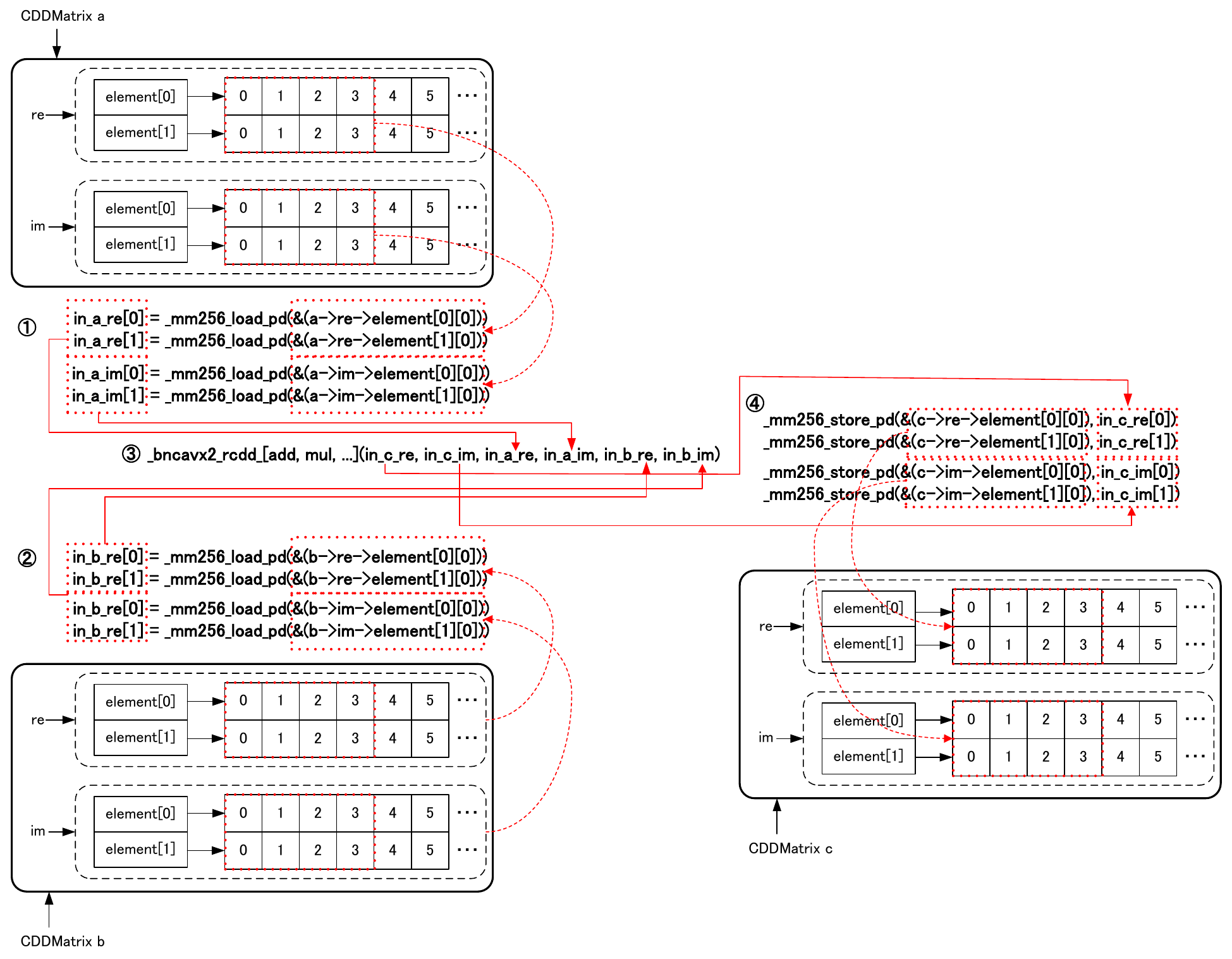}
    \caption{Complex linear computation with AVX2}\label{fig:complex_avx2}
    \end{center}
\end{figure*}

Therefore, SIMDization of complex arithmetic units is not necessary for basic matrix multiplication alone. However, in order to speed up complex LU decomposition, we decided to SIMDize the arithmetic, similar to the real LU decomposition, to speed up the element-by-element forward elimination calculation. For multi-component complex vectors and matrices, since the real and imaginary parts are separated, SIMD Load/Store instructions can be used to retrieve and store the elements. Using this approach, complex DD operations using AVX2 (\verb|_bncavx2_rcdd_add|, \verb|sub|, \verb|mul|, \verb|div|, \verb|inv|(reciprocal)) can be implemented as shown in \figurename\ref{fig:complex_avx2}. 

Our SIMDized complex arithmetic is based on a combination of SIMDized real arithmetic.

\subsection{Parallelization with OpenMP}

\begin{figure}[htb]
    \begin{center}
    \includegraphics[width=.6\textwidth]{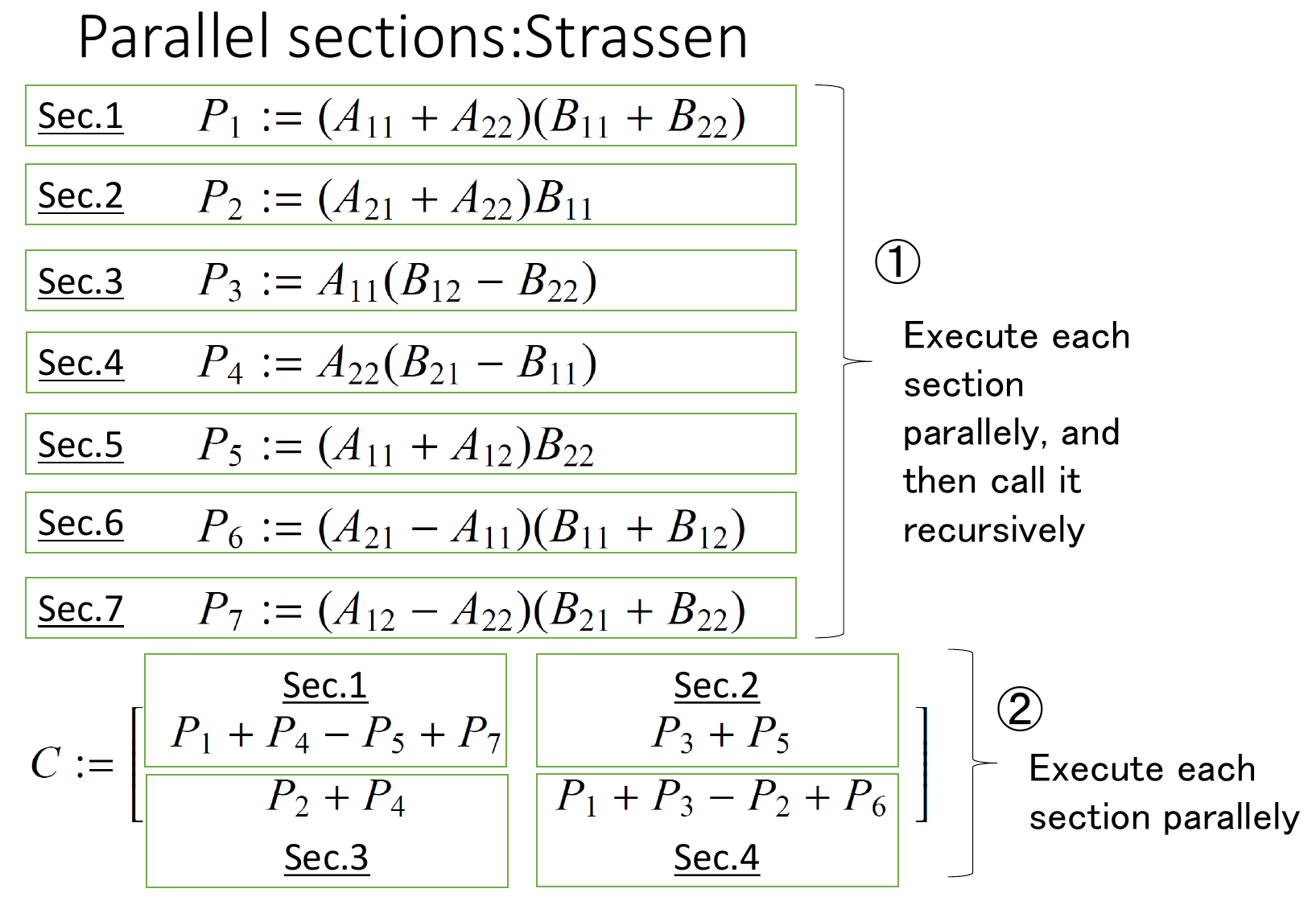}
    \caption{The parallelized Strassen algorithm}\label{fig:parallelized_strassen}
    \end{center}
\end{figure}

For parallelization, multi-component-type multiple precision complex linear computation is based on real linear computation units. The real matrix multiplication based on the Strassen algorithm has been parallelized in sections, as illustrated in \figurename\ref{fig:parallelized_strassen}. The parallelization is conducted section by section as depicted in. Due to the nested nature of the parallelization, performance improvement beyond 8 threads is not currently expected. Arbitrary-precision complex matrix multiplication is implemented using MPC operations, including parallelization.

\begin{figure}[htb]
    \begin{center}
    \includegraphics[width=.8\textwidth]{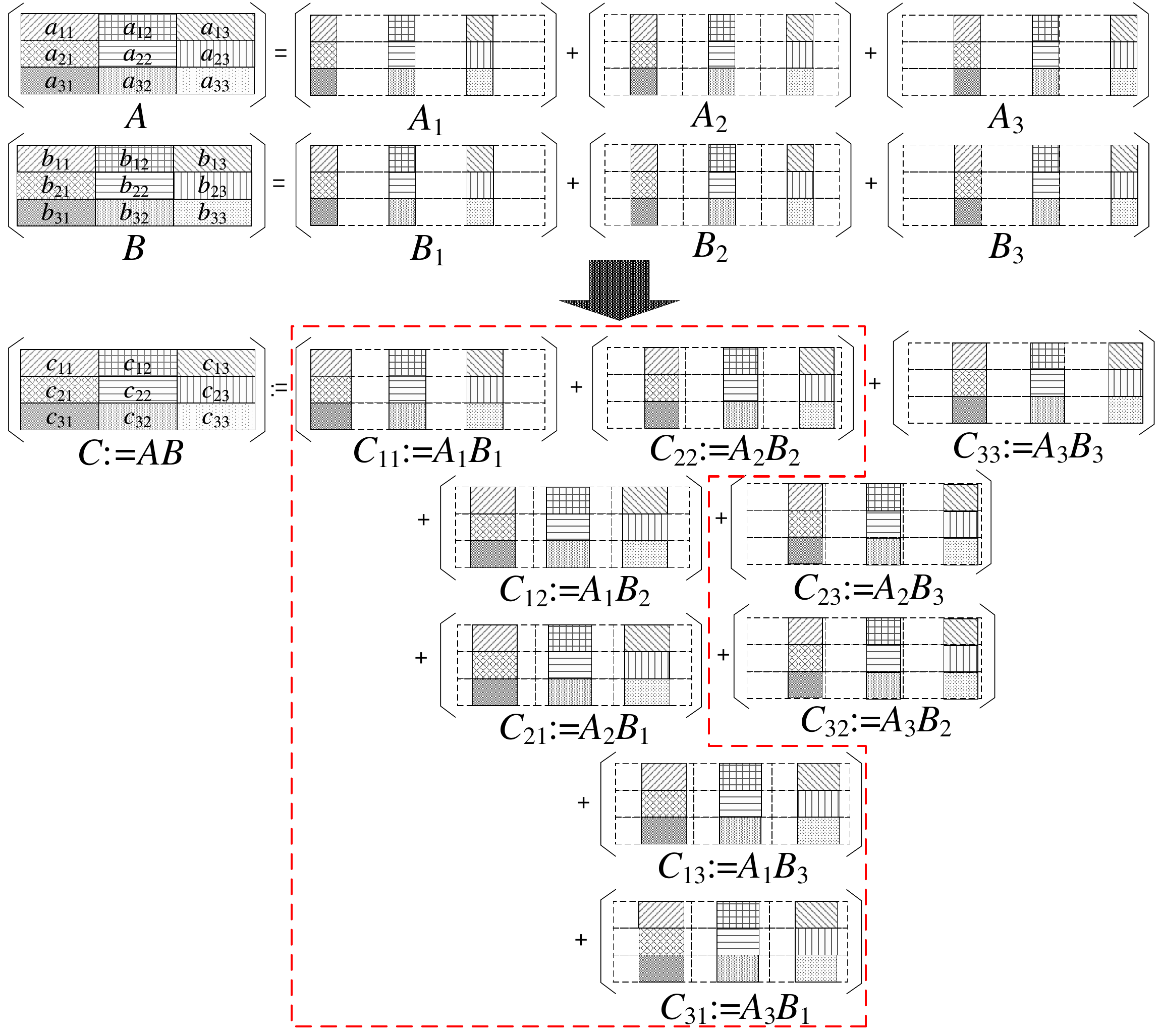}
    \caption{The parallelized Ozaki scheme}\label{fig:parallelized_ozaki_scheme}
    \end{center}
\end{figure}

The Ozaki scheme divides the computation into lengths to avoid introducing rounding errors in the base short-precision matrix multiplication (xGEMM, using binary64 DGEMM in this case). The product obtained by xGEMM is then added with a multiple-precision matrix addition to obtain a high-precision matrix product. To parallelize this process, it is necessary to use parallel xGEMMs or perform OpenMP loop parallelization for each xGEMM, as illustrated in \figurename\ref{fig:parallelized_ozaki_scheme}. In this case, the latter approach was used for parallelization.

\subsection{Benchmark test for complex matrix multiplication}

We use the following EPYC computing environment to present the results obtained through benchmark tests.
\begin{description}
    \setlength{\parskip}{0cm} 
    \item[CPU and OS] AMD EPYC 9354P 3.7 GHz 32 cores, Ubuntu 20.04.6 LTS
    \item[C/C++] Intel Compiler version 2021.10.0
    \item[Compiler Option] \verb|-O3 -std=c++11 -fp-model precise|
    \item[with AVX2] {\tt -axCORE-AVX2 -march=skylake -mtune=skylake -mcpu=skylake}
    \item[MPLAPACK] 2.0.1, GNU MP 6.2.1, MPFR 4.1.0, MPC 1.2.1
\end{description}

We address the problem of computing the matrix product $C := AB$ for given complex square matrices $A$ and $B$$\in\mathbb{C}^{n\times n}$. The elements of $A$ and $B$ used for benchmarking purposes are uniformly distributed random numbers in the interval $[0, 1]$. For the Ozaki scheme, we used a maximum of 6 divisions for DD accuracy, 8 divisions for TD accuracy, and 12 divisions for QD accuracy to ensure sufficient accuracy.

\begin{figure}[ht]
\begin{center}
    \includegraphics[width=.495\textwidth]{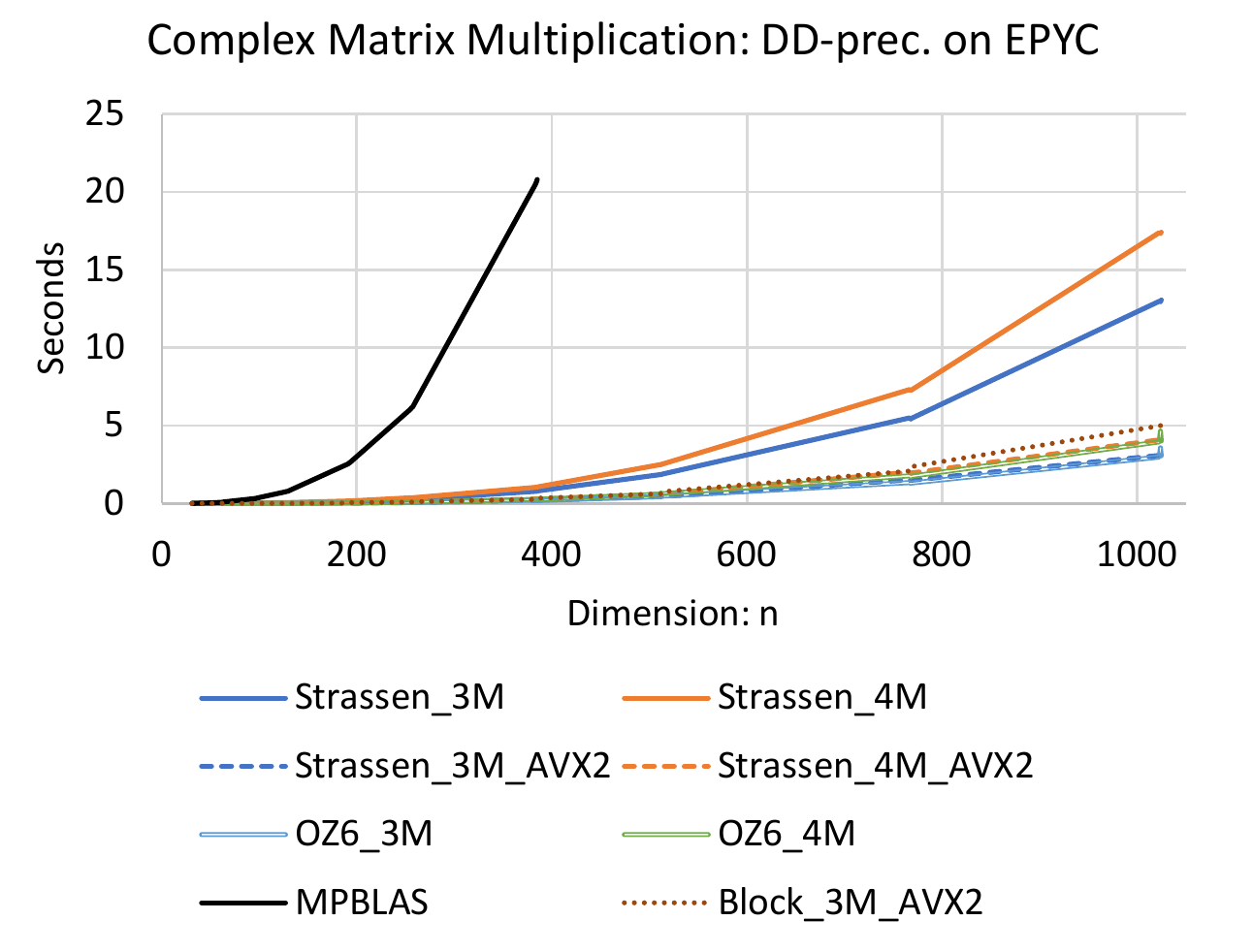}
    \includegraphics[width=.495\textwidth]{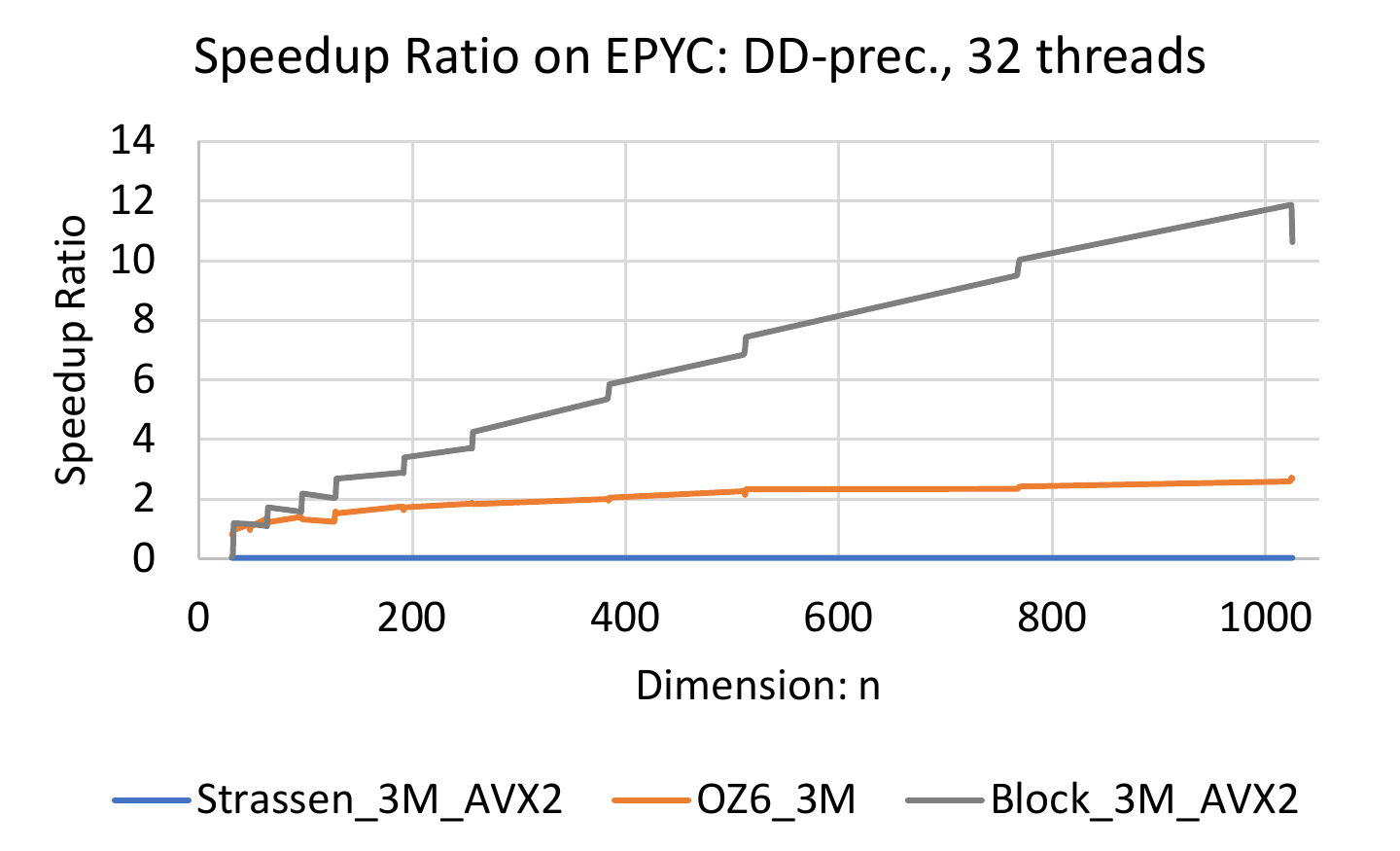}
    \caption{DD-prec.: Computation time in seconds (left) and speed-up ratio when using 32 threads (right)}\label{fig:cdd}
\end{center}
\end{figure}

At first, we show the results of DD-precision matrix multiplication on \figurename\ref{fig:cdd}. The left figure shows the computation time in seconds, and the right figure shows the corresponding speed-up ratio achieved by parallelization with OpenMP.

The left figure shows that our implementations of the Strassen and Ozaki schemes achieve up to 100 times performance improvement compared to the Cgemm routine of MPBLAS for serial calculations. Additionally, the 3M method exhibits a speed improvement of about 1.3 times over the 4M method, while the Strassen method achieves a speed-up of 3 to 4 times. The computation time of the Ozaki scheme is nearly equivalent to that of the AVX2 Strassen scheme due to the small number of partitions. Furthermore, the block matrix multiplication based on the 3M method is slightly slower than the Strassen matrix multiplication.

In OpenMP parallelization, the Strassen matrix multiplication does not exhibit any acceleration, even when using 32 threads. Similarly,  the Ozaki scheme demonstrates only a modest speed-up of about twice as fast. In contrast, the block matrix multiplication achieves a significant improvement, with speeds up to 12 times faster.

\begin{figure}[ht]
    \begin{center}
        \includegraphics[width=.495\textwidth]{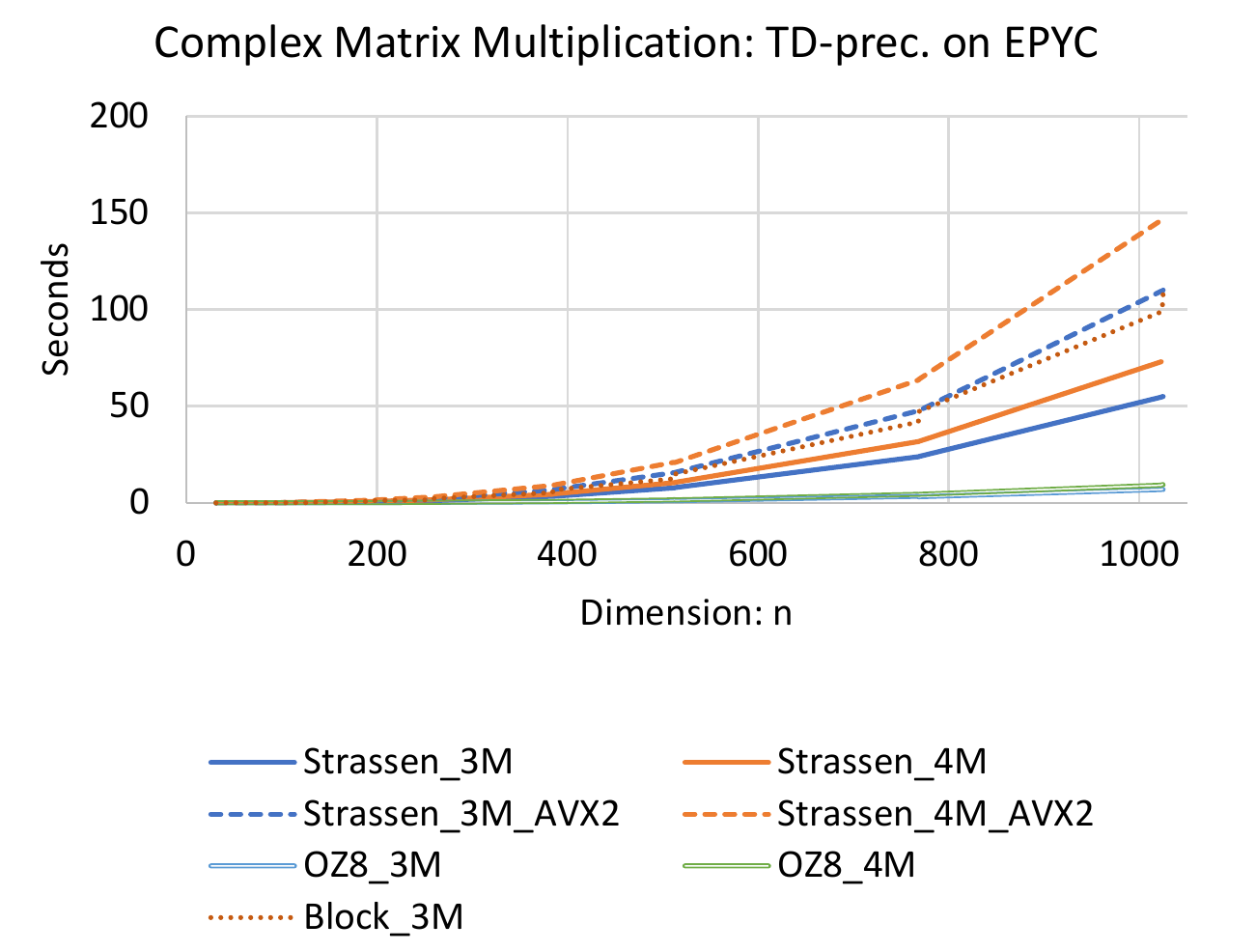}
        \includegraphics[width=.495\textwidth]{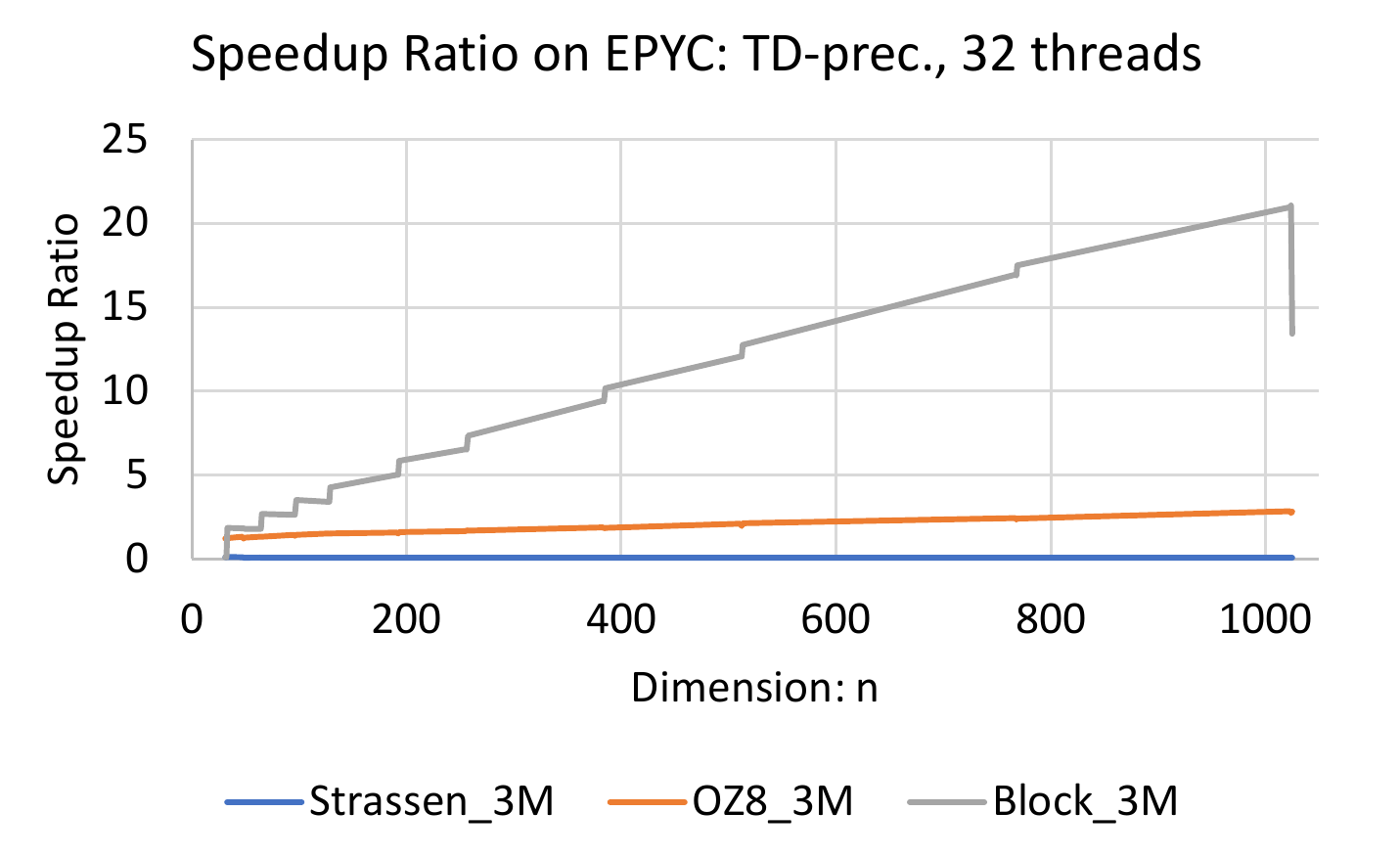}
        \caption{TD-prec. : Computation time in second (left) and speed-up ratio when using 32 threads (right)}\label{fig:ctd}
    \end{center}
\end{figure}

Next, we show the results obtained by TD-precision matrix multiplication in \figurename\ref{fig:ctd}.

The AVX2 SIMDized TD-precision Strassen, block matrix multiplication is observed to be slower than the non-SIMDized version, although the reason for this discrepancy is unclear. The 3M method demonstrates a speed improvement of about 1.3 times faster compared to the 4M method, which is common with DD precision. The Ozaki scheme achieves speed up to about 8 times faster than the Strassen 3M method.

In OpenMP parallelization, the Strassen matrix multiplication is observed to be slower, and the Ozaki scheme is only about 2 times faster, even when utilizing 32 threads, as in the DD accuracy. The block matrix multiplication is up to 20 times faster.

\begin{figure}[ht]
    \begin{center}
        \includegraphics[width=.495\textwidth]{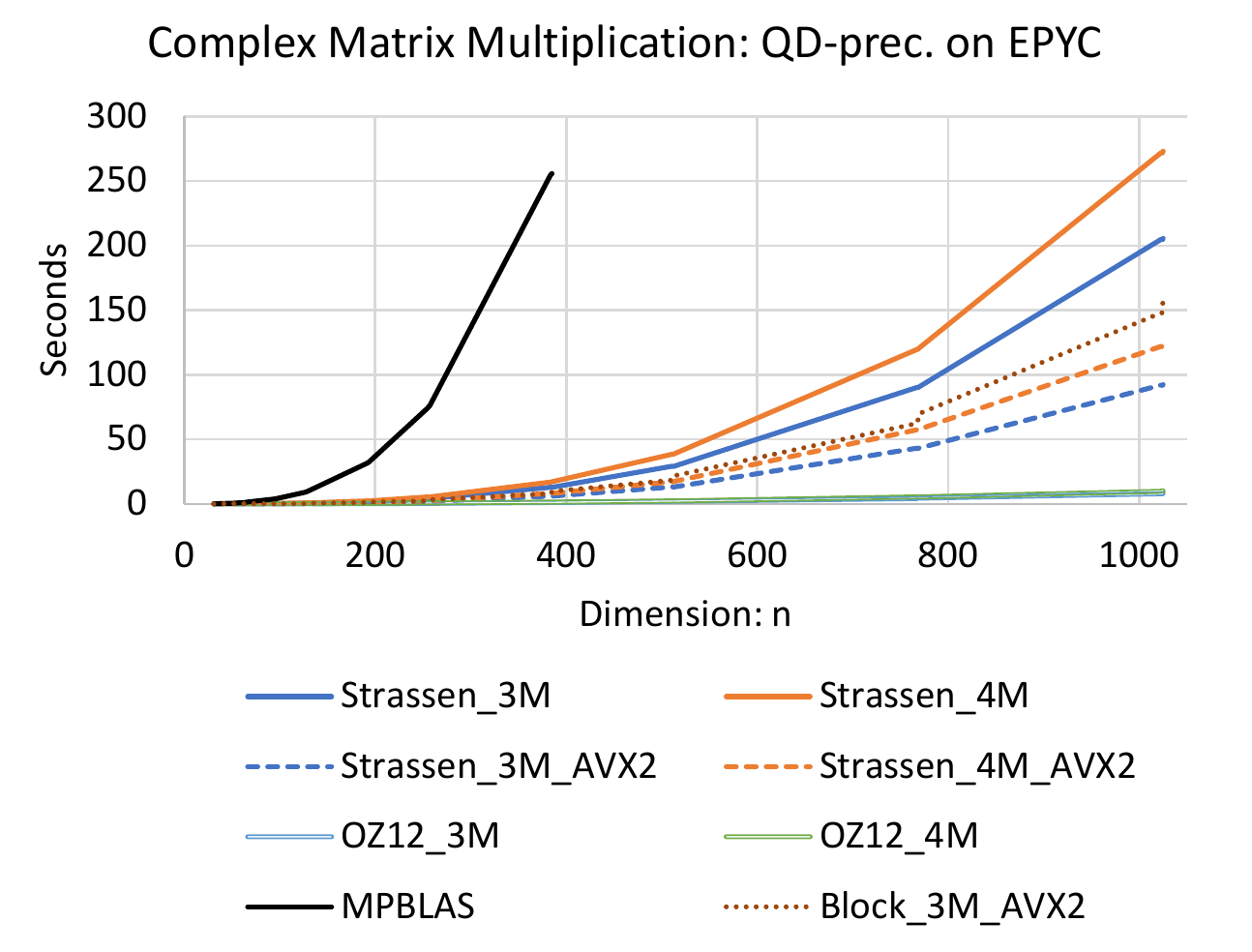}
        \includegraphics[width=.495\textwidth]{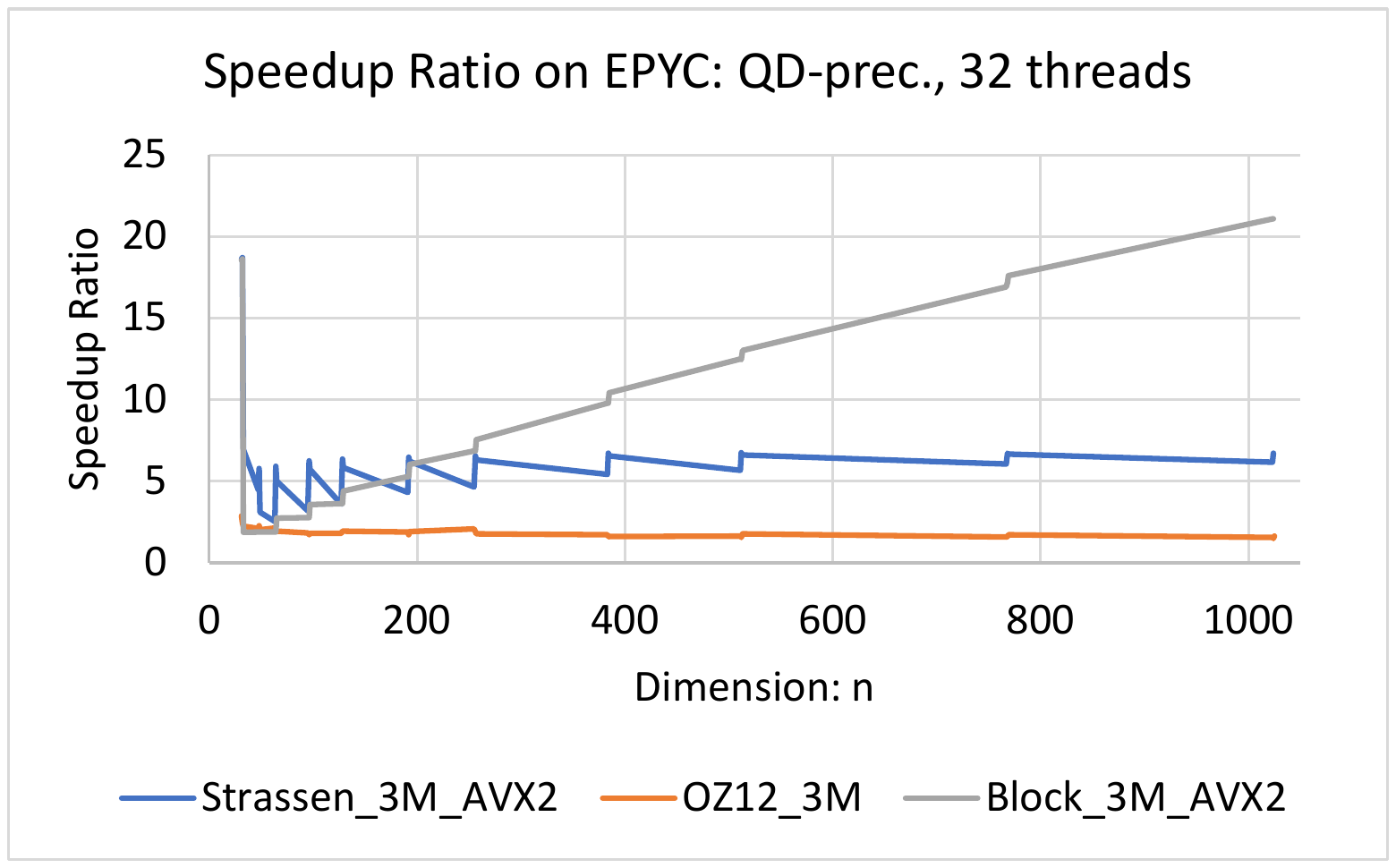}
        \caption{QD-prec. : Computation time in seconds (left) and speed-up ratio when using 32 threads (right)}\label{fig:cqd}
    \end{center}
\end{figure}

Finally, the results of QD-precision matrix multiplication are shown in \figurename\ref{fig:cqd}.

Similar to DD and TD precision, the Strassen matrix multiplication achieves approximately 2 times faster performance with AVX2 optimization and 1.3 times faster performance with the use of the 3M method.

In OpenMP parallelization, the Strassen matrix multiplication demonstrates faster performance, while the Ozaki scheme achieves only 2 times speed-up compared to the DD- and TD-precision schemes, despite utilizing 32 threads. The block matrix multiplication is 20 times faster than TD precision.

These results demonstrate that while the 3M method can speed up serial complex matrix multiplication, it does not speed up the parallelization in relation to the number of threads, except in the case of block matrix multiplication.
    
\section{Acceleration of multiple-precision complex LU decomposition}\label{sec:lu}

The complex LU decomposition is implemented based on a complex basic linear computation that integrates all the implemented speed-up techniques (3M method, SIMD, OpenMP parallelization) to achieve high performance. The following is a concise description of the implementation.

The target $n$-dimensional complex linear system is shown as 
\begin{equation}
	A\mathbf{x} = \mathbf{b}, \label{eqn:linear_eq}
\end{equation}
where $A\in\mathbb{C}^{n\times n}$, $\mathbf{x}\in\mathbb{C}^n$ and $\mathbf{b}\in\mathbb{C}^n$. We solve (\ref{eqn:linear_eq}) to determine the unknown vector $\mathbf{x}$ given $A$ and $b$.

\begin{figure}[htb]
    \begin{center}
        \includegraphics[width=.495\textwidth]{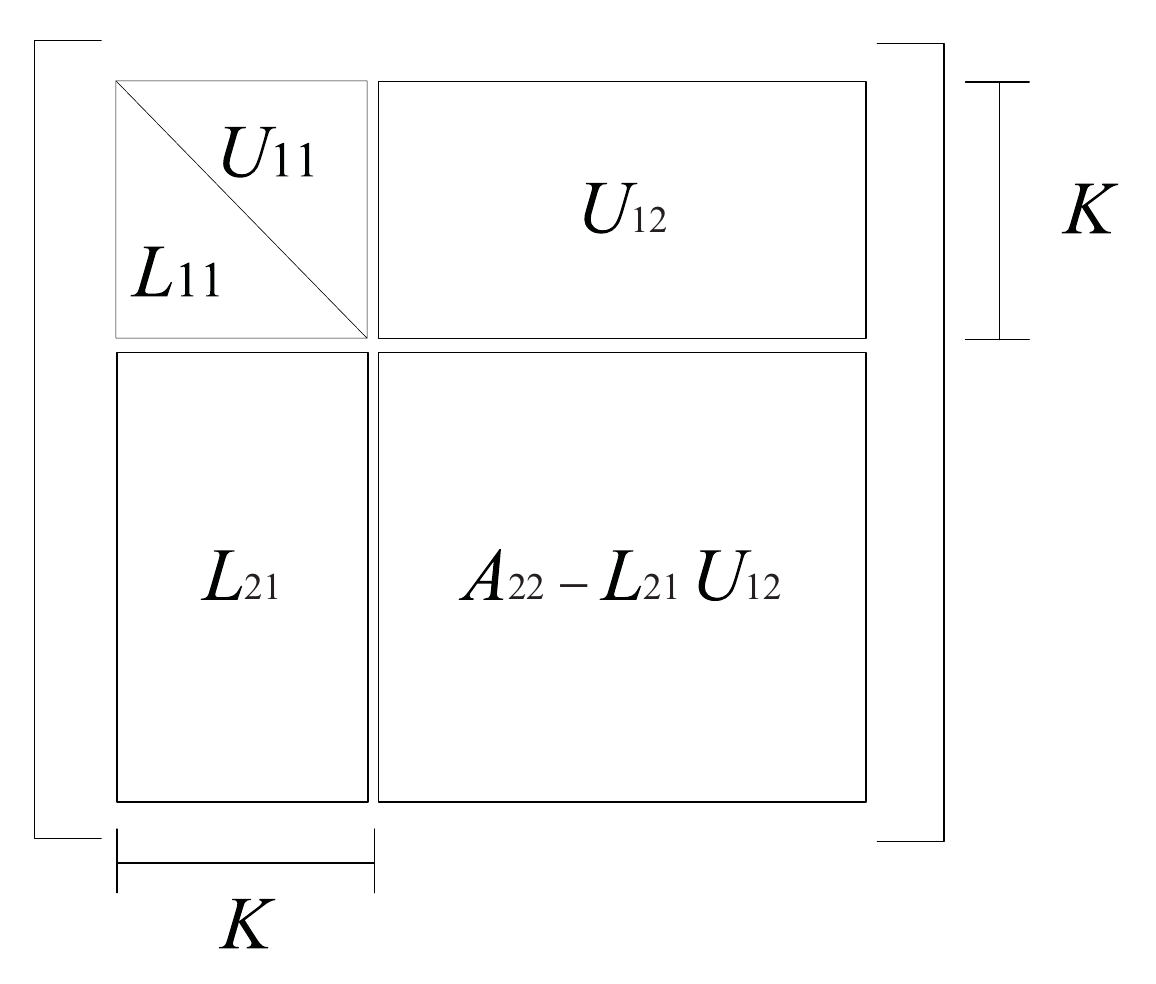}
        \caption{LU decomposition with matrix multiplication}
        \label{fig:lu}
    \end{center}
\end{figure}

The current standard implementation of LU decomposition uses fast matrix multiplication. As shown in \figurename\ref{fig:lu}, a certain block size $K$($=\mbox{NDIM}$) is set in advance, and $A - L_{12}U_{21}$ is performed to compute the $A_{22}$ part. Consequently, the computation performance varies depending on the value of $K$. The LU decomposition (Cgetrf function) of MPLAPACK used for comparison also  uses complex matrix multiplication (Cgemm function). In all cases, LU decomposition is computed using a partial pivoting strategy.

In addition, a standard LU decomposition approach (referred to as ``Normal LU" or ``Normal"), which eliminates rows sequentially without updating or performing matrix multiplication, was implemented using complex arithmetic with AVX2, as illustrated in \figurename\ref{fig:complex_avx2}, owing to the ease of speeding up the process with SIMD instructions. AVX2 demonstrated effectiveness in real LU decomposition, with the Ozaki scheme emerging as the fastest in this case. The question remains: will the same effect be obtained for complex LU decomposition ?

\section{Benchmark test of complex LU decomposition}

During the benchmark test, the computation time and maximum error of the numerical solution of the LU decomposition were measured by varying $K$ while calculating DD, TD, and QD precision, respectively. The graphs below show the computation time for normal LU and LU decomposition using matrix multiplication, as well as the relative error of the numerical solutions obtained through forward and backward substitutions.

Our focus is on solving the $n$-dimensional complex linear system depicted in equation (\ref{eqn:linear_eq}). For benchmarking purposes, the elements of $A$ are generated as uniform random numbers within the range $[0, 1]$. 

The true solution $\mathbf{x}$ is $(\mathbf{x})_k = k + k\cdot\mathrm{i}$, and the constant vector is derived by computing $\mathbf{b} := A\mathbf{x}$ using 2048-bit precision calculation.

\subsection{DD-precision computation}

First, the results of the DD accuracy calculation are presented. This involved measuring the relative error and computation time while increasing $1\leq K\leq n$ by 32.

\begin{figure*}[htb]
    \begin{center}
        \includegraphics[width=.495\textwidth]{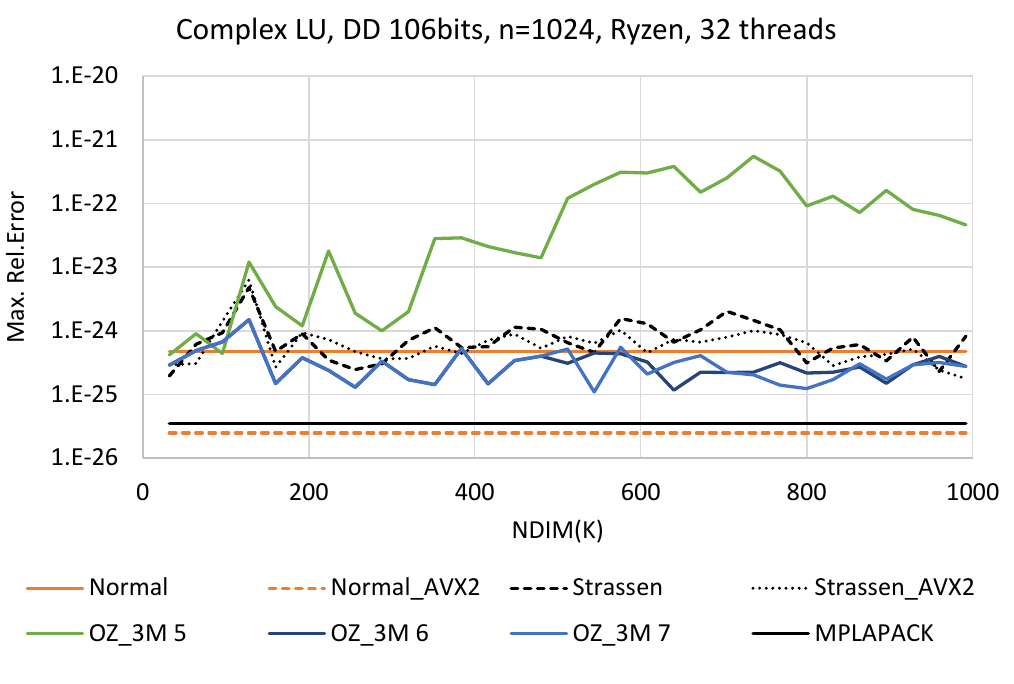}
    \caption{DD-precision complex LU decomposition: Max. relative error}\label{fig:clu_dd_relerr}
    \end{center}
\end{figure*}

Examining the relative error (\figurename\ref{fig:clu_dd_relerr}), the MPLAPACK and AVX2 SIMDized normal LU decompositions demonstrate superior performance. The normal LU decomposition and LU decomposition using the Ozaki scheme and Strassen matrix multiplication yield similar results. Overall, the results also show that the accuracy tends to decrease when the Strassen is used. For this problem, the Ozaki scheme with DGEMM does not achieve optimal accuracy with 5 divisions, necessitating more than 6 divisions. Therefore, results are shown for 5, 6, and 7 divisions.

\begin{figure*}[htb]
    \begin{center}
        \includegraphics[width=.495\textwidth]{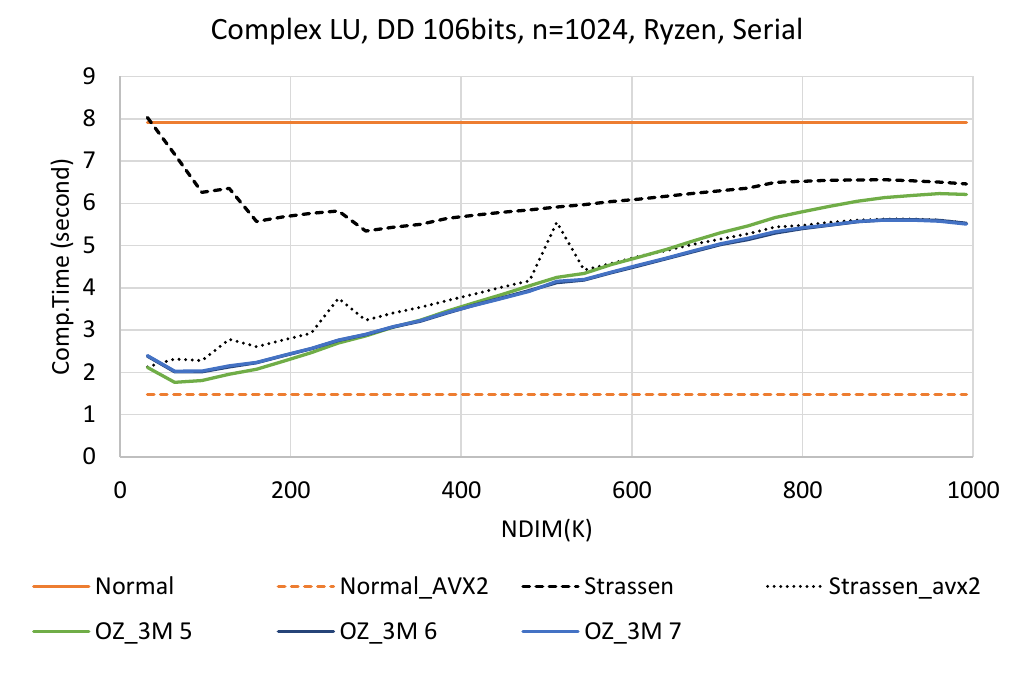}
        \includegraphics[width=.495\textwidth]{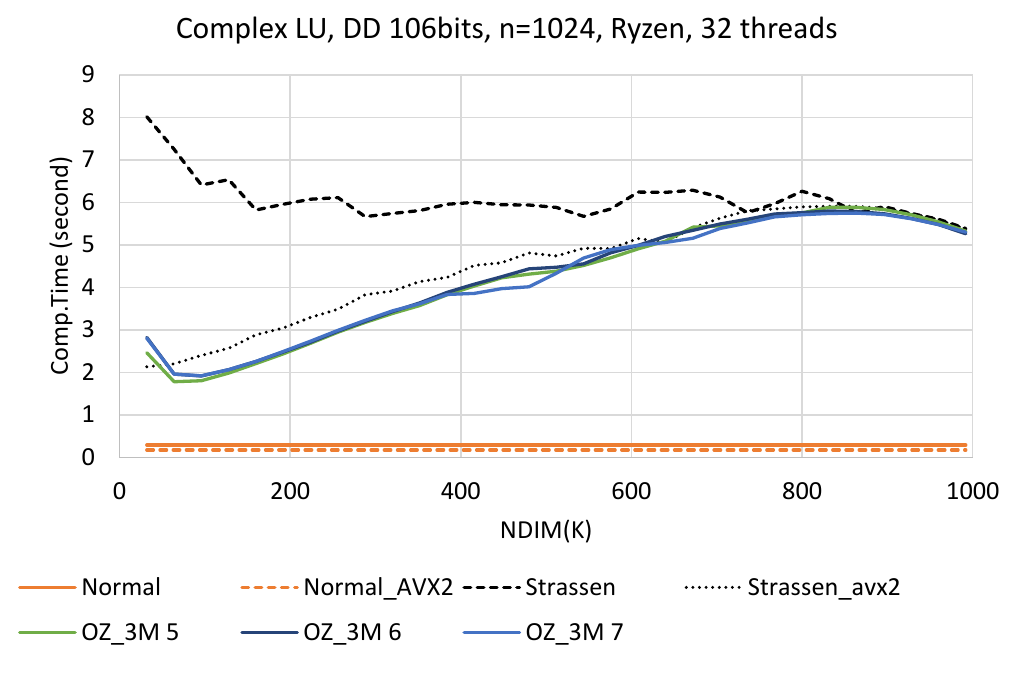}
    \caption{DD-precision complex LU decomposition: Seconds of serial computation (left), Seconds of 32-thread computation (right)}\label{fig:clu_dd_comptime}
    \end{center}
\end{figure*}

\figurename\ref{fig:clu_dd_comptime} demonstrates the measured serial computation time (left) and 32-thread parallel computation time (right). The results with normal LU decomposition and AVX2 acceleration are represented by solid orange and dashed lines, respectively. The relative error of MPLAPACK is depicted as a solid black line, although it is not included in the computation time graph due to the significant difference in scale (130.7 s).

The AVX2 normal LU decomposition exhibits the fastest serial computation, while the AVX2 Strassen LU decomposition and Ozaki scheme show similar computation times, with smaller $K$ values requiring less time. However, they are not as fast as the AVX2 normal LU decomposition.

\begin{figure}[htb]
    \begin{center}
        \includegraphics[width=.75\textwidth]{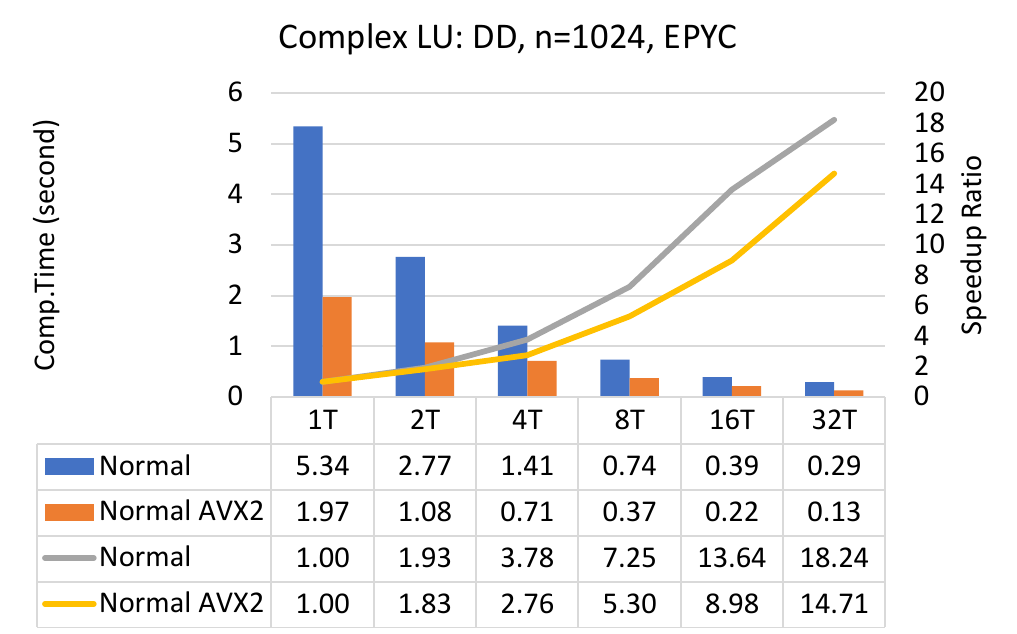}
    \caption{Parallelization of normal complex LU decomposition: DD-precision}\label{fig:normal_clu_dd_performance}
    \end{center}
\end{figure}

In the case of DD-precision calculations, OpenMP parallelization of matrix multiplication initially increases the computation time with 2 threads, after which it becomes faster proportionally to the number of threads. Therefore, with 32 threads, the computation time is approximately the same as for serial computation, and as a result, the parallelization effect is not enhanced. On the other hand, as depicted in \figurename\ref{fig:normal_clu_dd_performance}, the normal LU decomposition achieves speed-up in proportion to the number of threads, and the AVX2 version remains the fastest even with 32 threads.

The reason for this phenomenon is believed to be the higher computation complexity inherent in complex LU decomposition compared to real LU decomposition, despite the enhanced performance of complex operations achieved through SIMD. As previously mentioned, while SIMD algorithms are not essential for complex basic linear computation implementations, they may facilitate acceleration in algorithms such as the Fast Fourier Transform, which are required locally as LU decomposition.

\subsection{QD-precision computing}

Next, we present the results of the QD accuracy calculations.

\begin{figure*}[htb]
    \begin{center}\figurename\ref{fig:clu_qd_relerr} 
        \includegraphics[width=.495\textwidth]{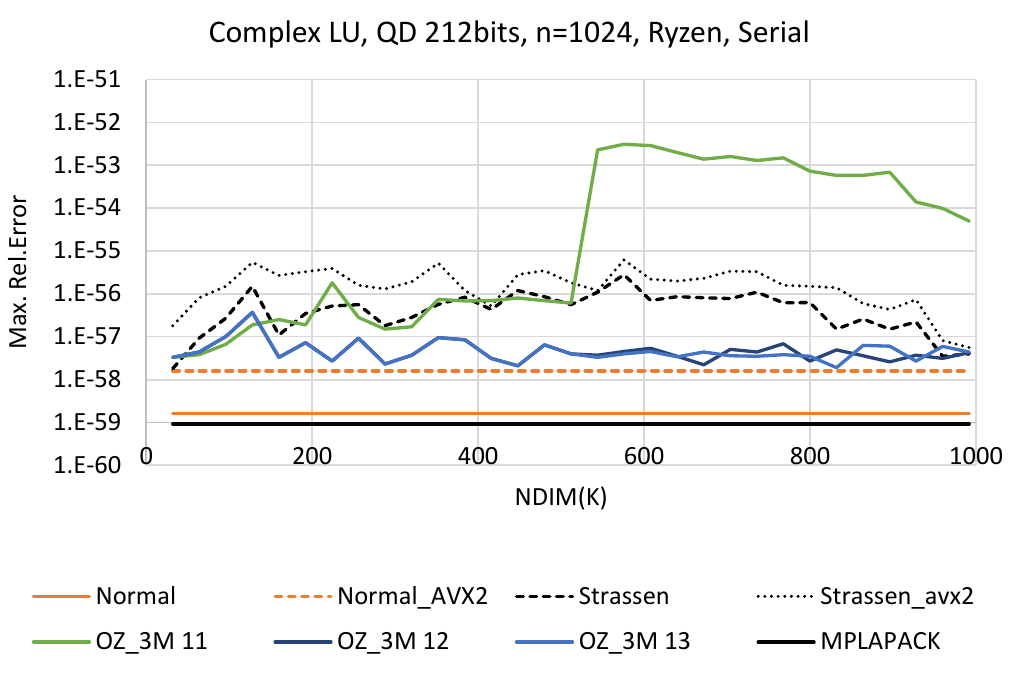}
    \caption{QD-prec. complex LU decomposition: Max. relative error}\label{fig:clu_qd_relerr}
    \end{center}
\end{figure*}

\figurename\ref{fig:clu_qd_relerr} displays the relative errors obtained from QD-precision benchmark tests. As depicted in this figure, the AVX2 SIMDized normal LU decomposition outperforms MPLAPACK, while the LU decomposition using the Ozaki scheme and the Strassen matrix multiplication exhibits larger relative errors compared to the normal LU decomposition. 

\begin{figure*}[htb]
    \begin{center}
        \includegraphics[width=.495\textwidth]{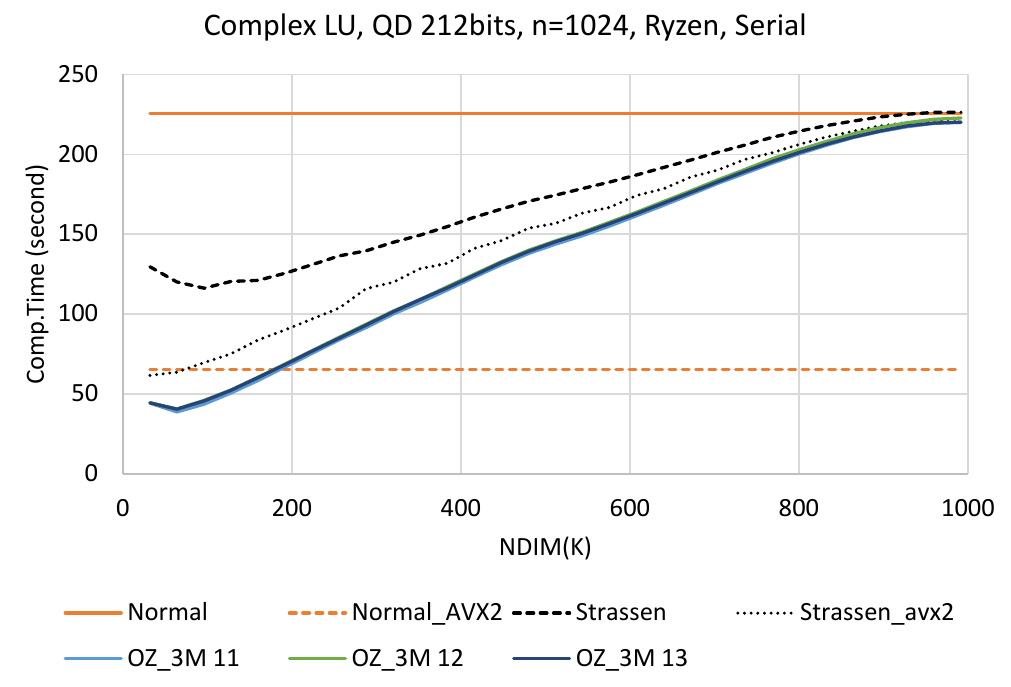}
        \includegraphics[width=.495\textwidth]{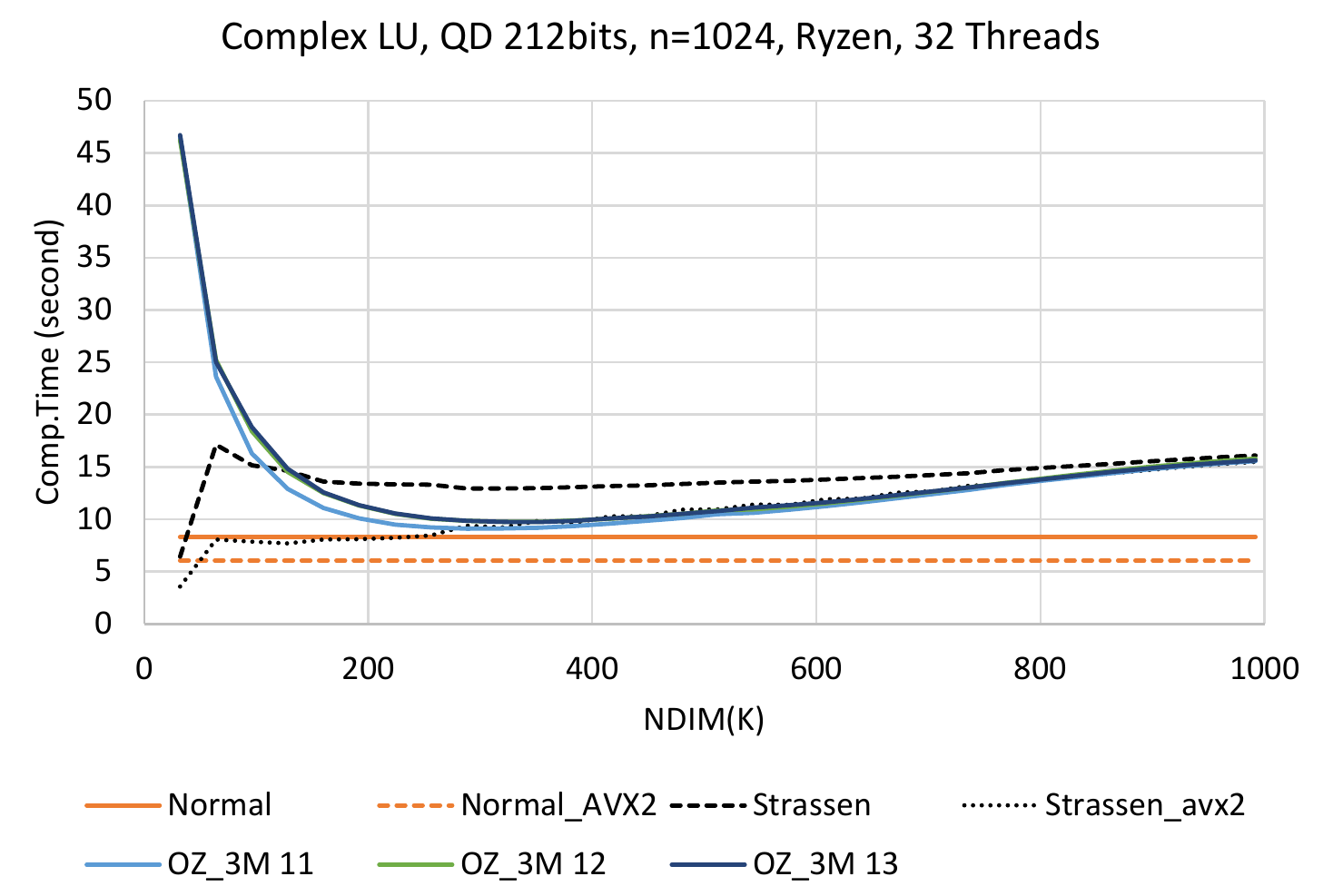}
    \caption{QD-prec. complex LU decomposition:Seconds of serial computation (left), Seconds of 32-thread computation (right)}\label{fig:clu_qd_comptime}
    \end{center}
\end{figure*}

\figurename\ref{fig:clu_qd_comptime} illustrates the computation time obtained from QD-precision benchmark tests. The graph showcases the measured serial computation time (left) and 32-thread parallel computation time (right) with the standard LU decomposition and the AVX2 accelerated approach, represented by solid and dashed orange lines, respectively. For this particular problem, the Ozaki scheme using DGEMM does not reach the highest QD accuracy with 11 divisions and requires more than 12 divisions, hence results for 11, 12, and 13 divisions are displayed. The relative error for MPLAPACK is represented by the solid black line, but the serial computation time is 595.4 seconds, so it is not included in the graph of computation time because it is on a different scale.

For serial calculations, the AVX2 SIMDized normal LU decomposition generally exhibits the fastest performance. However, the Ozaki scheme LU decomposition is slightly faster when $K$ is small. The Strassen LU decomposition is not as fast as the AVX2 SIMDized normal LU decomposition, even when AVX2 SIMDized.

\begin{figure}[htb]
    \begin{center}
        \includegraphics[width=.75\textwidth]{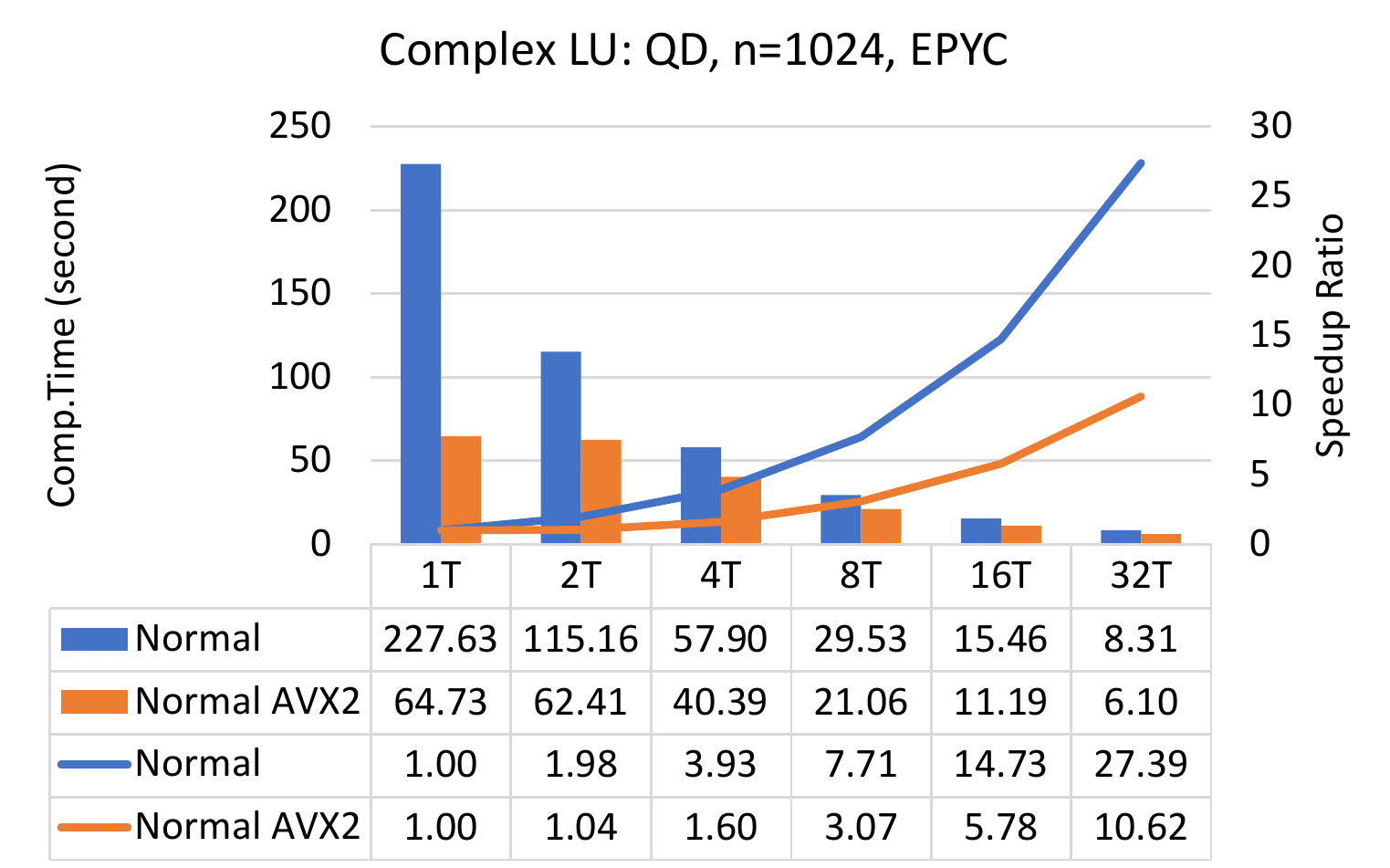}
    \caption{Parallelization of normal complex LU decomposition: QD-prec.}\label{fig:normal_clu_qd_performance}
    \end{center}
\end{figure}

In QD-precision calculations, the OpenMP parallelization of matrix multiplication accelerates the computation linearly with the number of threads. However, the parallelization effect of the normal LU decomposition surpasses this, with the AVX2 SIMDized normal LU decomposition remaining the fastest overall, except in certain cases where the AVX2 SIMDized Strassen LU decomposition outperforms it. As illustrated in \figurename\ref{fig:normal_clu_qd_performance}, although the parallelization efficiency is high, it is slowed down by AVX2.

\subsection{TD-precision computing}

Finally, we show the results of the TD-precision calculation, which represents an intermediate accuracy level between DD and QD.

\begin{figure*}[htb]
    \begin{center}
        \includegraphics[width=.495\textwidth]{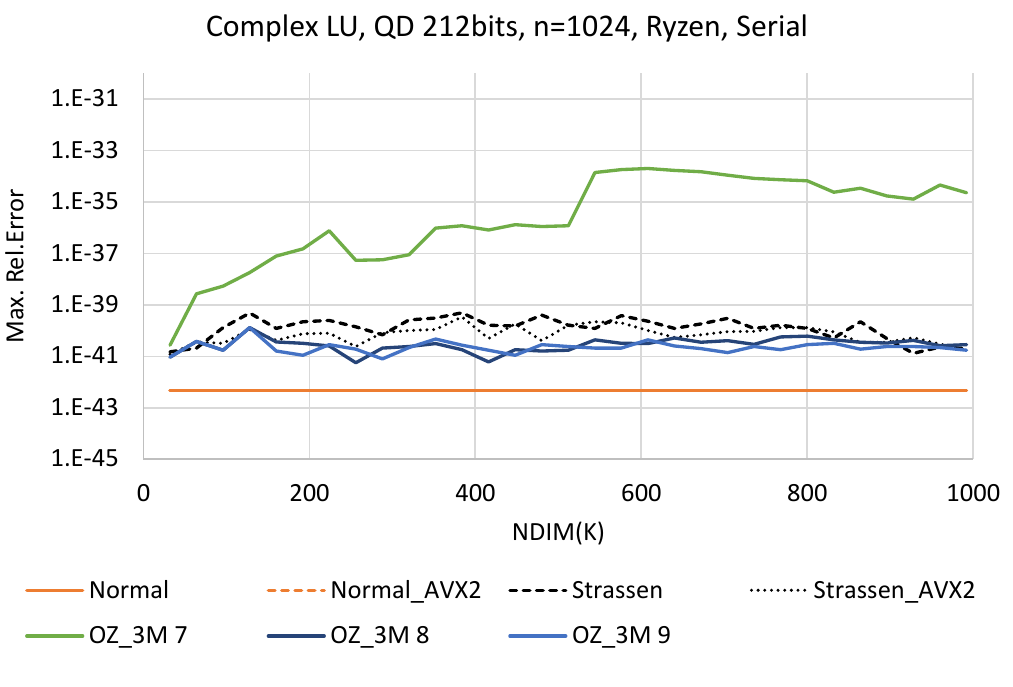}
        \caption{TD-precision complex LU decomposition: Max. relative error}\label{fig:clu_td_relerr}
    \end{center}
\end{figure*}

In this case, as depicted in \figurename\ref{fig:clu_td_relerr}, for TD accuracy, the Ozaki scheme with DGEMM does not achieve the highest accuracy with 7 divisions, requiring 8 or more.

The relative error analysis reveals that, similar to the case of QD accuracy, the AVX2-ized normal LU decomposition performs the best. Furthermore, the relative error observed in the LU decomposition using the Ozaki scheme and Strassen matrix multiplication is greater compared to the normal LU decomposition.

\begin{figure*}[htb]
    \begin{center}
        \includegraphics[width=.495\textwidth]{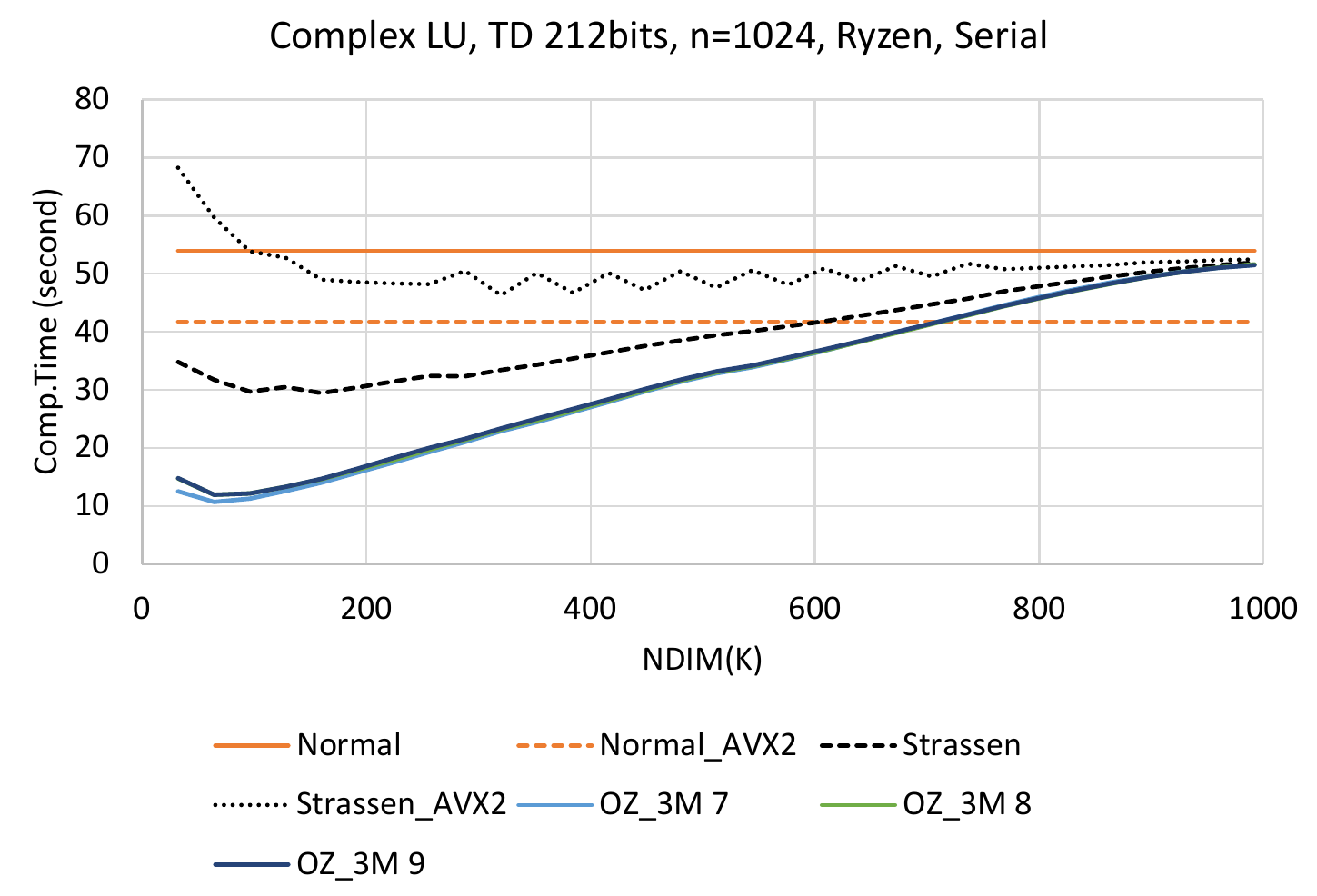}
        \includegraphics[width=.495\textwidth]{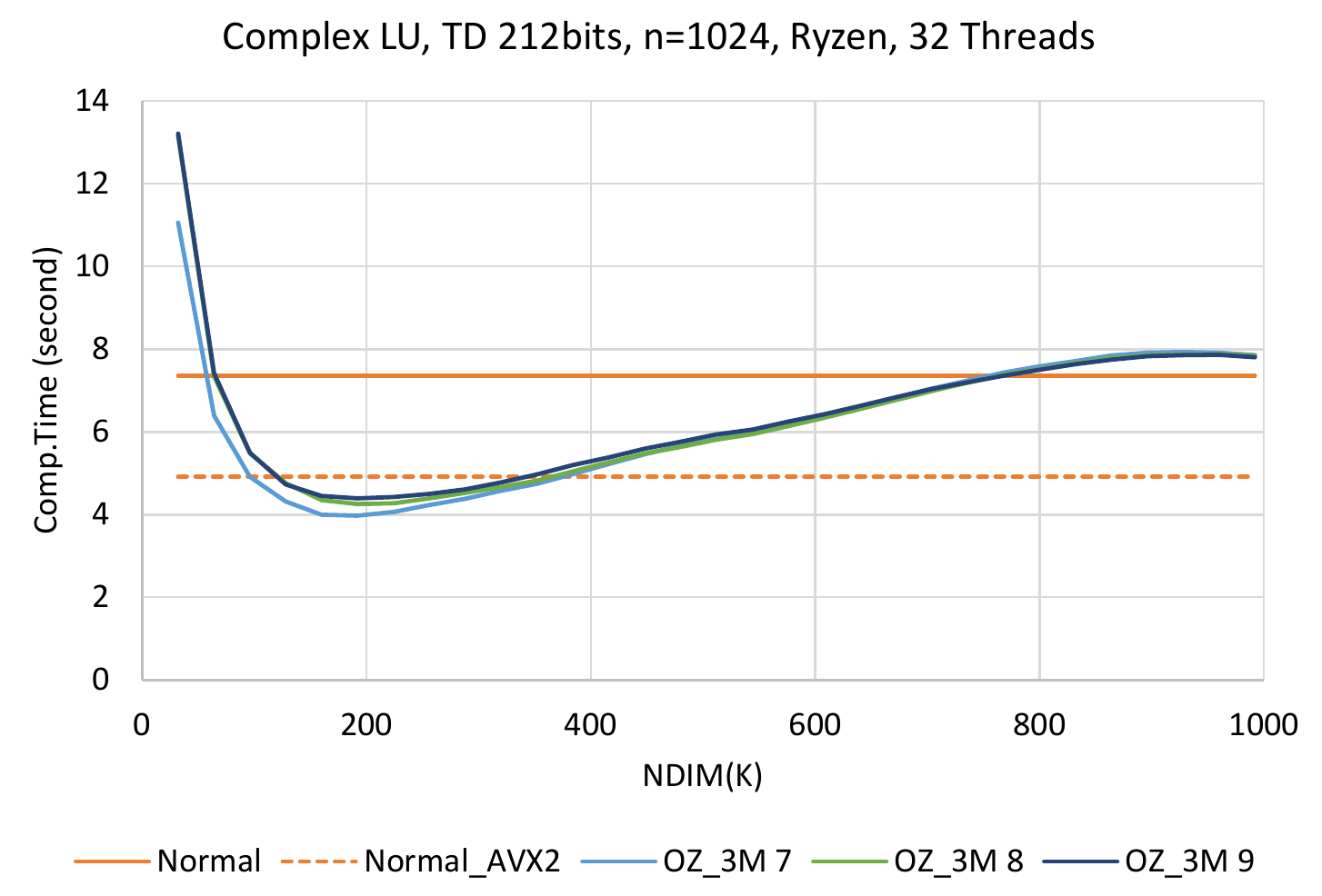}
    \caption{TD-precision complex LU decomposition: Seconds of serial computation (left), Seconds of 32-thread computation (right)}\label{fig:clu_td_comptime}
    \end{center}
\end{figure*}

A graph depicting computation time is presented in \figurename\ref{fig:clu_td_comptime}, omitting the 32-thread parallelization using TD accuracy due to its significantly increase in computation time for small $K$.

\begin{figure}[htb]
    \begin{center}
        \includegraphics[width=.75\textwidth]{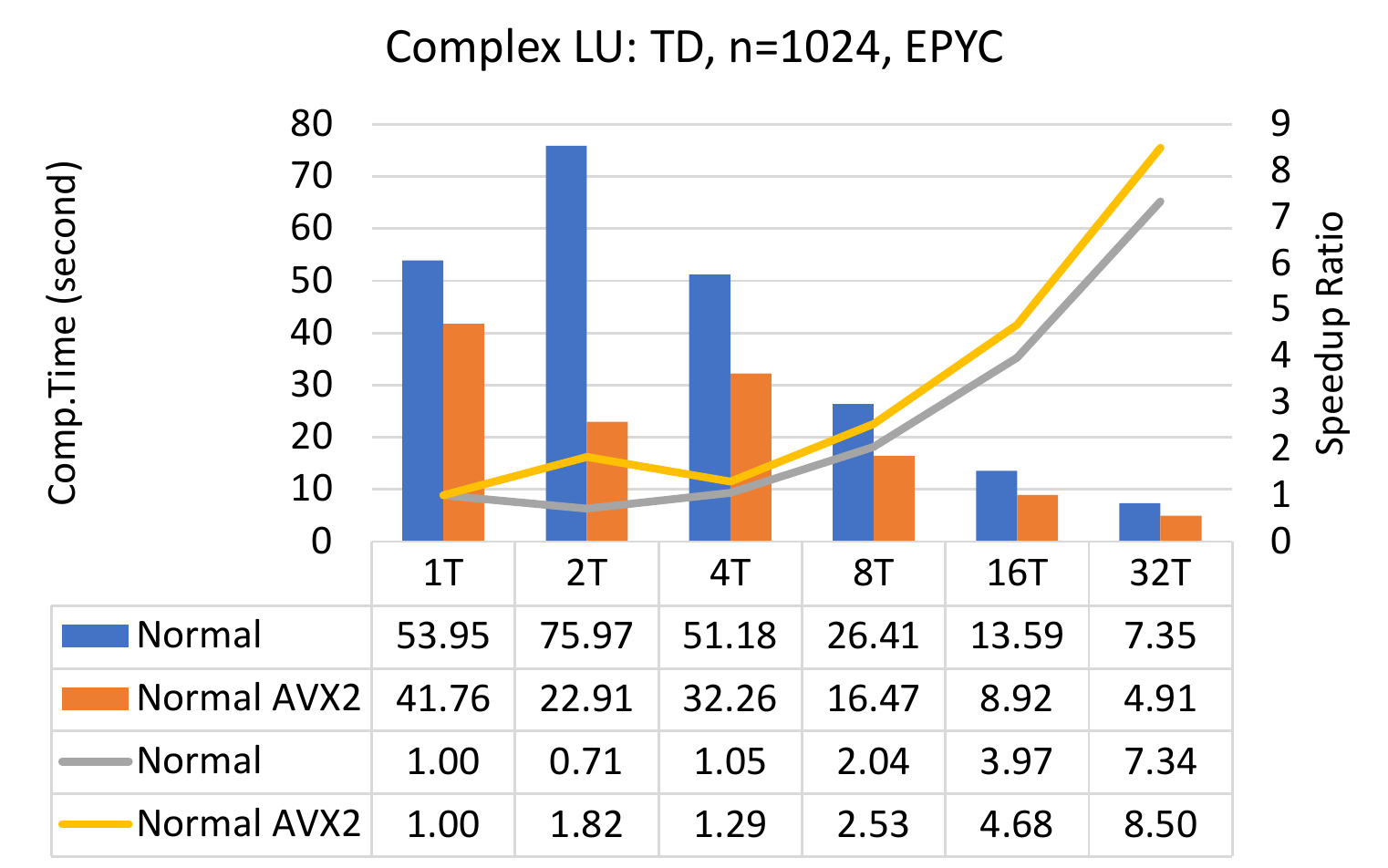}
    \caption{Parallelization of normal complex LU decomposition: TD-precision}\label{fig:normal_clu_td_performance}
    \end{center}
\end{figure}

For serial calculations, the AVX2 SIMDized normal LU decomposition generally does not outperform others, with the LU decomposition using the Ozaki scheme being largely faster for small $K$ values. In TD-precision calculations, OpenMP parallelization of the matrix multiplication displays straightforward speed-up proportionate to the number of threads. However, the parallelization effect of the normal LU decomposition outweighs it. The normal LU decomposition with AVX2 remains the fastest, except in cases where LU decomposition with the Ozaki scheme shows better performance. As shown in \figurename\ref{fig:normal_clu_td_performance}, the parallelization efficiency is considerable; however, the ratio is slowed down compared with DD and QD normal LU decomposition.

\subsection{Summary of high-performance complex LU decomposition}

These benchmark results indicate that the AVX2 SIMDized normal LU decomposition generally exhibits the fastest performance in terms of both DD and QD accuracy, with good parallelization. A summary solely focusing on the computation time is presented in \tablename\ \ref{table:comptime_clu}, highlighting the minimum computation time achieved for Strassen with AVX2 for DD and QD, and without AVX2 in TD.

\begin{table}[htb]
    \begin{center}
        \caption{Minimum computation time of complex LU decomposition in second}\label{table:comptime_clu}
        \begin{tabular}{|c|c|c|c|c|c|}\hline
        \multicolumn{6}{|c|}{DD precision, $n=1024$} \\ \hline
                        & MPLAPACK  & Normal    & Normal AVX2 & Strassen AVX2 & OZ \\ \hline
        Serial          & 130.7     & 7.9       & 1.5         & 2.1 ($K=32$)  & 2.0 ($K=64$)\\    
        Parallel(32T)   & $-$       & 0.29      & 0.18        & 2.1 ($K=32$)  & 1.9 ($K=96$) \\ \hline\hline
        \multicolumn{6}{|c|}{TD precision, $n=1024$} \\ \hline
        Serial          & $-$       & 54.0      & 41.8        & 46.7 ($K=384$) & 11.9($K=64$)\\    
        Parallel(32T)   & $-$       & 7.4       & 4.9         & 8.5 ($K=992$)  & 4.3($K=192$)\\ \hline\hline
        \multicolumn{6}{|c|}{QD precision, $n=1024$} \\ \hline
        Serial          & 595.4     & 227.6     & 64.7        & 60.1($K=32$)   & 40.7($K=64$)\\    
        Parallel(32T)   & $-$       & 8.3       & 6.1         & 3.6($K=32$)    & 9.8($K=352$)\\ \hline
        \end{tabular}
    \end{center}
\end{table}

Finally, it was demonstrated that the 32-thread AVX2 normal LU decomposition outperformed MPLAPACK with serial computation by approximately 726 times for DD accuracy, and approximately 91 times for QD accuracy. These findings suggest that the default setting for complex LU decomposition should be the AVX2 SIMDized normal LU decomposition.

\section{Conclusion and future works}

The implementation of the multiple precision complex basic linear computation, based on the enhanced multiple-precision real basic linear computation, demonstrated improved speed by using the 3M method for matrix multiplication. Additionally, the complex LU decomposition was found to be sufficiently fast through the utilization of SIMD complex operations and parallelization with OpenMP. However, it was observed that the speed-up of the LU decomposition is significant enough to surpass conventional LU decomposition algorithms in DD- and QD-precision, given the current complex matrix multiplication, both in serial and parallel computations.

Future work will focus on enhancing the performance of numerical algorithms for ill-conditioned nonlinear problems by improving the following functions. 

\begin{enumerate}
    \item Implementation of a multiple precision sparse matrix and extension of our implementation to use real and complex sparse matrices.
    
    \item Speed up reproducible computations for binary64 and other types of precision arithmetic.
    
    \item Extending fast, multi-precision computation to the Python environment. 
\end{enumerate}
    
This includes seeking further speed-up of real and complex basic linear computation of dense matrices and LU decomposition. Of particular interest is confirming the potential speed-up of SIMD-enhanced multi-component-type real and complex arithmetic for other numerical algorithms.

%
\subsubsection{Acknowledgment}
This research is supported by the JSPS Grant-in-Aid for Scientific Research 23K11127. 

%
%
%
%

\end{document}